\documentclass[11pt,a4paper]{article}
 \usepackage{euscript}

\usepackage{amsmath}
\usepackage{graphicx}
\usepackage{amsthm}
\usepackage{amssymb}
\usepackage{epsfig}
\usepackage[all]{xy}
\usepackage{url}

\theoremstyle{theorem}
\newtheorem{theorem}{Theorem}[section]

\newtheorem{lemma}[theorem]{Lemma}
\newtheorem{proposition}[theorem]{Proposition}

\newtheorem{remark}{Remark}[section]
\numberwithin{equation}{section}

\title{On additive and multiplicative arithmetical functions}
\date{\today}
\author{Masood Aryapoor \footnote{masood.aryapoor@mail.ipm.ir}
\footnote{This research was in part supported by a grant from IPM.}
\\ \\
School of Mathematics, \\
Institute for Research in Fundamental Sciences (IPM),\\
P.O.Box: 19395-5746, Tehran, Iran.
}

\begin{document}
\maketitle
\begin{abstract}
\noindent
A characterization of multiplicative (and additive) arithmetical functions is given.  Using 
this characterization, we show that the group of multiplicative arithmetical functions is 
isomorphic to the group of additive arithmetical functions. 
\end{abstract} 
\begin{section}{Introduction}
It is known that the set of multiplicative arithmetical functions is a group under the (Dirichlet) product. 
It is clear that the set of additive arithmetical functions is a group under addition. The main purpose of 
this paper is to show that these two groups are isomorphic. In order to show this, we use 
some operations defined on   arithmetical functions, especially those defined in \cite{CS}.
\end{section}
\begin{section}{Arithmetical functions}
Here is a review on arithmetical functions, for details see \cite{Ap} for example. An arithmetical function
is simply a function $a:\mathbb{N}\to \mathbb{C}$. The sum and the product of two arithmetical functions $a,b$ are defined via
$$(a+b)(n)=a(n)+b(n),$$
$$(a\star b)(n)=\sum_{d|n}a(d)b(n/d).$$
The set of arithmetical functions is a commutative $\mathbb{C}$-algebra under the above operations
(the scalar multiplication is defined by $(ra)(n)=ra(n)$ for an arithmetical function $a$ and  $r\in\mathbb{C}$).
The identity of this algebra is denoted by $I$ which is  the arithmetical function $I(1)=1$, $I(n)=1$ for $n>1$.
Given an arithmetical function $a$, we denote the product $a\star a\star ...\star a$ (with $n$ terms) by
$a^{*n}$.
Some arithmetical functions are given below:\\
(1) $u(n)=1$ for every $n\in\mathbb{N}$.\\
(2) $\mu=u^{-1}$: the Mobius function.\\
(3) $\phi$: the Euler totient function.\\
(4) $\Lambda$: the Mangoldt function.\\
(5) $\lambda$:   Liouville's function  function.\\
(6) $\sigma_c$: divisor functions, where $\sigma_0$ is denoted by $d$.\\
(7) $\nu(n)=$ the number of distinct prime divisors of $n$.\\
(8) $N(n)=n$  for every $n\in\mathbb{N}$.\\
(9) $\Omega(n)=$ the total number of prime divisors of $n$. \\
The derivative $a'$ of an arithmetical function $a$ is defined by $a'(n)=a(n)\log(n)$.
An arithmetical function $a$ is called multiplicative if $a(1)=1$ and $a(mn)=a(m)a(n)$ for every relatively prime natural numbers $m,n$. An arithmetical function $a$ is called
completely multiplicative if $a(1)=1$ and $a(mn)=a(m)a(n)$ for every natural
numbers $m,n$.  An arithmetical function $a$ is called additive if $a(mn)=a(m)+a(n)$
for every relatively prime natural numbers $m,n$. An arithmetical function $a$ is
called completely additive if $a(mn)=a(m)+a(n)$ for every natural numbers $m,n$.

One can define various operations on arithmetical functions, see \cite{Cs,Re} for example. We consider two of them. One of them is $\log:I+\mathbb{M}\to \mathbb{M}$ defined by
$$\log(a)=\sum_{n=1}^{\infty}\frac{(-1)^{n-1}(a-I)^{*n}}{n}.$$
Another one is $\exp:\mathbb{M}\to I+\mathbb{M}$ defined by
$$\exp(a)=I+\sum_{n=1}^{\infty}\frac{a^{*n}}{n!}.$$
It is known that both these operations are bijective and $\log=\exp^{-1}$, see \cite{CS}. Moreover
$\log(a\star b)=\log(a)+\log(b)$ and $\exp(a+b)=\exp(a)\star\exp(b)$. In other words,
$\log$ (and $\exp$) are group isomorphisms (here $\mathbb{M}$ is a group with addition and
$I+\mathbb{M}$ is a group with multiplication). It is also easy to see that
$\log(a)'=a'\star a^{-1}$ and $\exp(a)'=a'\star \exp(a)$.

\end{section}

\begin{section}{A representation of arithmetical functions}
Let $p_1=2,p_2=3,p_3,...$ be the sequence of prime numbers. 
One can see that there is  a unique   valuation $v$ (with the multiplicative group $\mathbb{Q}^{+}$ as its value group) on $\mathbb{C}(x_1,x_2,...)$
such that $v(x_1^{\alpha_1}\cdots x_n^{\alpha_n})=p_1^{\alpha_1}\cdots p_n^{\alpha_n}$ for an arbitrary  monomial
$ x_1^{\alpha_1}\cdots x_n^{\alpha_n}\in \mathbb{C}(x_1,x_2,...)$. The completion of $\mathbb{C}(x_1,x_2,...)$ with respect to this valuation is denoted by
$\mathbb{K}$ and its corresponding valuation ring is denoted by $\mathbb{O}$. In order to present elements of $\mathbb{K}$, we introduce some notations. Let $r$ be a positive rational number. Then $r$ can be uniquely
written as $p_1^{\alpha_1}\cdots p_n^{\alpha_n}$ for some integers $\alpha_i$'s.
Set $[\textbf{x}]^r:=x_1^{\alpha_1} x_2^{\alpha_2}\cdots\in \mathbb{C}(x_1,x_2,...)$.  Note that
$v([\textbf{x}]^r)=r$ and $[\textbf{x}]^r [\textbf{x}]^s=[\textbf{x}]^{rs}$.
Any element
$f\in\mathbb{K}$ can be uniquely written as
$$f(x_1,x_2,...)=\sum_{r\in \mathbb{Q}^{+}} a_r [\textbf{x}]^r,$$
(where $a_r\in\mathbb{C}$) such that there is some $r_0\in \mathbb{Q}^{+}$ depending on $f$ with
$\{r| a_r\neq 0\}\subset r_0\mathbb{N}$.  We call $\{r| a_r\neq 0\}$ the support of $f$ and denote it by
$\text{Supp}(f)$. Since the support of $f$ is a well-ordered subset of $\mathbb{Q}$ we have
$\text{Supp}(f)=\{r_1,r_2,r_3,...\}$ where $r_1<r_2<r_3<...$ are positive rational numbers. So $f$ can be written as
$$f=\sum_{i} a_{r_i} [\textbf{x}]^{r_i}.$$
It is easy to see that $f\in\mathbb{O}$ iff   $\text{Supp}(f)\subset \mathbb{Q}^{\geq 1}$.
Let $\mathbb{A}$ be the set of elements $f$ of $\mathbb{K}$ such that $\text{Supp}(f)\subset \mathbb{N}$.
It is clear that $\mathbb{A}$ is a ring and its field of fractions is just $\mathbb{K}$. Any element of $f\in \mathbb{A}$ can be uniquely written
as $$f=a_1+a_2[\textbf{x}]^{2}+\cdots+a_n[\textbf{x}]^{n}+\cdots.$$
Clearly $\mathbb{A}$ is the completion of $\mathbb{C}[x_1,x_2,...]$  with respect to $v$. Moreover $\mathbb{A}$ is a local ring with maximal ideal
$\mathbb{M}$ consisting of those $f\in \mathbb{A}$ with $a_1=0$.
The following proposition is well known, see \cite{CE}.
\begin{proposition}
The map $a\mapsto a(1)+a(2)[\textbf{x}]^{2}+\cdots+a(n)[\textbf{x}]^{n}+\cdots$ gives an isomorphism of
rings between the ring of arithmetical functions and $\mathbb{A}$.
\end{proposition}
We identify the ring of arithmetical functions and $\mathbb{A}$ using the above isomorphism from now on.
More precisely given an arithmetical function $a:\mathbb{N}\to \mathbb{C}$, we denote $ a(1)+a(2)[\textbf{x}]^{2}+\cdots+a(n)[\textbf{x}]^{n}+\cdots$ by $a(x_1,x_2,...)$.

The following identities are easy to check (the limits are taken in the topology of $\mathbb{K}$ defined by $v$)
$$u(x_1,x_2,...)=\prod_{i=1}^{\infty}\frac{1}{1-x_i},\quad \mu(x_1,x_2,...)=\prod_{i=1}^{\infty}(1-x_i),$$ $$\phi(x_1,x_2,...)=\prod_{i=1}^{\infty}\frac{1-x_i}{1-p_ix_i},\quad
\lambda(x_1,x_2,...)=\prod_{i=1}^{\infty}\frac{1}{1+x_i},$$
$$\Lambda(x_1,x_2,...)=\sum_{i=1}^{\infty}\frac{\log(p_i)x_i}{1-x_i},\quad d(x_1,x_2,...)=\prod_{i=1}^{\infty}\frac{1}{(1-x_i)^2},$$
$$N(x_1,x_2,...)=\prod_{i=1}^{\infty}\frac{1}{1-p_ix_i},\quad \sigma_{c}(x_1,x_2,...)=\prod_{i=1}^{\infty}\frac{1}{(1-x_i)(1-p_i^{c}x_i)},$$
$$\nu(x_1,x_2,...)=(\prod_{i=1}^{\infty}\frac{1}{1-x_i})\sum_{i=1}^{\infty}x_i,\quad \Omega(x_1,x_2,...)=(\prod_{i=1}^{\infty}\frac{1}{1-x_i})\sum_{i=1}^{\infty}\frac{x_i}{1-x_i}.$$
In the following lemma, a characterization of multiplicative (and additive) arithmetical functions is given.
\begin{lemma}\label{Lcharacterization}
(1) An arithmetical   function $a$ is multiplicative iff there are (unique) formal power series $f_i(x)\in\mathbb{C}[[x]]$ whose constant terms are 1 such that $a(x_1,x_2,...)=\prod_{i=1}^{\infty}f_i(x_i)$.
An arithmetical   function $a$ is completely multiplicative iff
there are complex numbers $c_1,c_2,...$ such that
$a(x_1,x_2,...)=\prod_{i=1}^{\infty}(1-c_ix)^{-1}$, i.e. $a=u(c_1x_1,c_2x_2,...)$.\\
(2) An arithmetical   function $a$ is additive iff there are (unique) formal power series $f_i(x)\in\mathbb{C}[[x]]$ whose constant terms are 0 such that
$$a(x_1,x_2,...)=u(x_1,x_2,...)\sum_{i=1}^{\infty}f_i(x_i).$$
An arithmetical   function $a$ is completely additive iff
there are complex numbers $c_1,c_2,...$ such that
$$a(x_1,x_2,...)=u(x_1,x_2,...)\sum_{i=1}^{\infty}\frac{c_ix_i}{1-x_i}.$$
\end{lemma}
\begin{proof}
(1) If $a$ is multiplicative then
$$a(x_1,x_2,...)=\prod_{i=1}^{\infty}(\sum_{n=0}^{\infty} a(p_i^n)x_i^n).$$ Conversely,
if $f_i(x)=\sum_{j=0}^{\infty} c_{ij}x^j$, then $a(p_1^{m_1}...p_k^{m_k})= \prod_{i=1}^{k}c_{im_i}$. Hence $a$
is multiplicative. The uniqueness of $f_i$'s follows from $a(p_i^n)=c_{in}$. Finally, $a$ is completely
multiplicative iff there are complex numbers  $c_1,c_2,...$ with $c_{ij}=c_i^j$. So $a$
is completely multiplicative iff   $f_i(x)=\sum_{j=0}^{\infty} (c_{i}x)^j$ iff $f_i(x)=(1-c_ix)^{-1}$.\\
(2) If $a$ is additive then
$$a(x_1,x_2,...)=u(x_1,x_2,...)\sum_{i=1}^{\infty}((1-x_i)\sum_{n=1}^{\infty} a(p_i^n)x_i^n).$$ So
if $a$ is additive then
$a(x_1,x_2,...)= u(x_1,x_2,...)\sum_{i=1}^{\infty}f_i(x_i)$ with $f_i(x)=a(p_i)x+(a(p_i^2)-a(p_i))x^2+\cdots$. Conversely, if
$f_i(x)(1-x)^{-1}=\sum_{j=1}^{\infty} c_{ij}x^j$ then $a(p_1^{m_1}...p_k^{m_k})= \sum_{i=1}^{k}c_{im_i}$.
The uniqueness of $f_i$'s follows from $a(p_i^n)=c_{in}$. Finally, $a$ is completely
additive iff there are complex numbers  $c_1,c_2,...$ with $c_{ij}=jc_i$. So $a$
is completely multiplicative iff   $f_i(x)(1-x)^{-1}=\sum_{j=0}^{\infty} jc_{i}x^j=c_ix(1-x)^{-2}$, i.e. $f_i(x)=c_ix(1-x)^{-1}$.
\end{proof}
This gives the following characterization of additive functions.
\begin{proposition}
An arithmetic function $a$ is additive iff $\mu\star a(n)=0$ for every natural number $n$ with $\nu(n)\neq 1$,
i.e. $n$ is not a power of a prime number.
\end{proposition}
\begin{proof}
If $a$ is additive then, by Lemma \ref{Lcharacterization},   there are (unique) formal power series $f_i(x)\in\mathbb{C}[[x]]$ whose constant terms are 0 such that
$$a(x_1,x_2,...)=u(x_1,x_2,...)\sum_{i=1}^{\infty}f_i(x_i).$$
so $(\mu\star a)(x_1,x_2,...)=\sum_{i=1}^{\infty}f_i(x_i)$. This clearly shows that $\mu\star a(n)=0$ for every natural number $n$ with $\nu(n)\neq 1$. Conversely, the property $\mu\star a(n)=0$ for every natural number $n$ with $\nu(n)\neq 1$ implies that
 $(\mu\star a)(x_1,x_2,...)=\sum_{i=1}^{\infty}f_i(x_i)$ for some formal power series $f_i(x)\in\mathbb{C}[[x]]$ whose constant terms are 0. So $$a(x_1,x_2,...)=u(x_1,x_2,...)\sum_{i=1}^{\infty}f_i(x_i)$$ is additive
 by Lemma \ref{Lcharacterization}.
\end{proof}

\end{section}

\begin{section}{Isomorphism of two groups}
The set $M$ of multiplicative arithmetical functions is an abelian group under multiplication.
Also, the set $A$ of additive arithmetical functions is an abelian group under addition. These two groups
are isomorphic as proved in the following theorem.  
\begin{theorem}
The map $\Psi:I+\mathbb{M}\to\mathbb{M}$ defined by
$\Psi(a)=u\star \log(a)$ is an isomorphism of groups. Moreover $a$ is multiplicative iff $\Psi(a)$ is additive.
Therefore $\Psi:M\to A$ is an isomorphism of groups.
\end{theorem}
\begin{proof}
Since multiplication by $u$ is a group isomorphism of $\mathbb{M}$ and 
$\log:I+\mathbb{M}\to\mathbb{M}$ is a group isomorphism (see \cite{CS}), we see that
$\Psi:I+\mathbb{M}\to\mathbb{M}$  is also a group isomorphism. So we only need to prove the last part.   Suppose that   $a$ is multiplicative. So, by Lemma \ref{Lcharacterization}, there are (unique) formal power series $f_i(x)\in\mathbb{C}[[x]]$ whose constant terms are 1 such that $a(x_1,x_2,...)=\prod_{i=1}^{\infty}f_i(x_i)$. Then
$$\Psi(a)(x_1,x_2,...)=u(x_1,x_2,...)\sum_i\log(f_{i}(x_{i})),$$
which is additive  by Lemma \ref{Lcharacterization}.   Conversely let $\Psi(a)$ be additive. Then, by Lemma \ref{Lcharacterization}, there are (unique) formal power series $f_i(x)\in\mathbb{C}[[x]]$ whose constant terms are 0 such that
$$\Psi(a)(x_1,x_2,...)=u(x_1,x_2,...)\sum_{i=1}^{\infty}f_i(x_i).$$
This implies that $\log(a)(x_1,x_2,...)= \sum_{i=1}^{\infty}f_i(x_i)$. So
$$a(x_1,x_2,...)=\exp(\log(a)(x_1,x_2,...))=\exp(\sum_{i=1}^{\infty}f_i(x_i))=
\prod_{i=1}^{\infty}\exp(f_i(x_i)),$$
is multiplicative, by Lemma \ref{Lcharacterization}.
\end{proof}
It is easy to see that the inverse of $\Psi$ is given by $\Psi^{-1}(a)=\exp(\mu\star a)$. 
\begin{remark}
In \cite{Re}, the author shows that some of the abelian groups ($M$ being one of them) defined using various classes of arithmetical functions are isomorphic. Using the above theorem,
we see that the group of additive arithmetical functions
is also isomorphic to those groups considered in \cite{Re}.
\end{remark}

\end{section}

\end{document}